\documentstyle{article}

\def \E {{\bf {E\,}}}

\def \d {\delta}

\def \v {\vert}

\def \la {\langle}

\def \ra {\rangle}

\def \T {{\hbox{ Tr}}}

\def \g {\gamma}

\def \o {\omega}

\def \O {\Omega}

\author{A.Khorunzhy\\
UFR de Math\'ematiques \\
Universit\'e Paris 7 - Denis Diderot, FRANCE\thanks{Permanent address:
 Mathematical Division, Institute for Low Temperature Physics, Kharkov,
UKRAINE, \quad e-mail: khorunjy@ilt.kharkov.ua}}

\title{PRODUCTS OF RANDOM MATRICES AND q-CATALAN NUMBERS
}

\begin{document}

\maketitle

\centerline{ }
\begin{abstract}

\noindent We give an interpretation of the q-Catalan numbers 
in frameworks of the random matrix theory and weighted 
partitions of the set of integers.

\end{abstract}

 
\hskip 1.9cm non-commutative random variables


\vskip 0.5cm

In a joint  discussion \cite{CM}, Christian Mazza asked, 
what one can obtain when regarding the weighted pairings
of
$2k$ points under the condition that they are non-crossing? In the present
note we give one possible answer to this question.

Let us
consider a set of $N$-dimensional random matrices
{\mbox{$A^{(r)}, \  r\in {\bf N}$}} determined on the same probability space.
We assume that these matrices are real symmetric and
$$
[A_{N}^{(r)}]_{ij} = {1\over \sqrt N} a_{ij}^{(r)}, \quad i,j = 1,\dots, N,
$$
where $\{a_{ij}^{(r)}, i\le j\}$ is the family of Gaussian 
random variables with zero mathematical expectation
and the covariance matrix
$$
\E a_{ij}^{(r)} a_{i'j'}^{(r')} = V(r-r')
\big( \d_{ii'}\d_{jj'} + \d_{ij'} \d_{i'j}
\big).
\eqno (1)
$$
Relations (1) mean that 
the probability distribution of the matrix 
$A^{(r)}$ with given $r$ is equal to that of GOE
\cite{M}. Certainly we assume $V$ to be real, even, and positively 
defined function.

In this note, our main subject is the product
$$
Q_{2k}^{(N)}(p) = \prod_{s=1}^{2k} A^{(s)}_N
$$
and we will be related mainly with the simplest case when
$$
V(r-r')= p^{\v r-r'\v}\ , \ 0<p<1.
\eqno (2)
$$

If $N=1$, we obtain the product of real gaussian random variables.
In this case we can use the integration by parts
formula for a Gaussian random vector $\vec \g=(\g_1, \dots, \g_n)$
with zero mean value:
$$
\E \g_l F(\vec \g) = \sum_{m=1}^n \E \g_l\g_m \E {\partial F(\vec \g)\over
\partial \g_m},
\eqno (3)
$$
where $F$ is a nonranom function.
Then the mathematical expectation
$$
\E Q^{(1)}_{2k} (p) = \sum_{\o \in \O_{2k}} \prod_{(l,m)\in \o} p^{\v l-m\v}
\eqno (4)
$$
is given by the sum of weighted partitions of the set $(1,2,\dots,2k)$ of 
$2k$ labelled points  into $k$ pairs. Here
$\o$ denotes one particular partition,
$\O_{2k}$ stands for the set of all possible partitions,
and $(l,m)\in \o$ means that given a partition $\o$,
the product is taken over pairs that form $\o$ and  whose elements are denoted
by
$(l,m)$.
 
Certainly, if $p=1$, then 
$$
\E Q^{(1)}_{2k}(1) = {(2k)!\over 2^k k!}.
$$
It is interesting  to study (4) in the limit $k\to\infty$,
namely,
the asymptotic behaviour of the value
$$
{1\over k} \ \log \left\{\E Q^{(N)}_{2k} (p)\right\}, \ \ N=1, \ k\to\infty
\eqno (5)
$$
in the dependence of the parameter $p$.
It is proved that there exists critical value $p_c$
such that (5) is positive for  $p>p_c$
and negative for $p<p_c$ \cite{CM}.

Now let us turn to the case of non-commutative random variables,
namely, to the limit  of $N=\infty$.

\vskip 1cm
{\bf Theorem 1.1}

{\it The average value of $Q_{2k}^{(N)}(p) $}
$$
  {1\over N} \T\left\{
A^{(1)}A^{(2)} \cdots A^{(2k)}\right\} = \la Q_{2k}^{(N)}(p) \ra  
\eqno (6)
$$
{\it converges as $N\to\infty$ in average
$$
\lim_{N\to\infty}\ \E \la Q_{2k}^{(N)}(p) \ra=B_k(p);
\eqno (7)
$$
the limit is  given by equality
$$
B_{k}(p)=  p\phi_k(p^2),
$$
where $\phi_k(p)$ are  determined by 
recurrent relations 
$$
\phi_{k+1}(p) = \sum_{i=1}^{k+1} \ p^i \, \phi_{i-1}(p) \, \phi_{k+1-i}(p),
\quad 
$$
with initial conditions $\phi_0(p)=1,\ \phi_1(p)=1$. The numbers
$C_k(q)= \phi_k(p)\ p^{ {{k-1}\choose{2}}}$, where  $q=p^{-1}$,  are known as
the q-Catalan numbers
\cite{S}.}

\vskip 1cm
To discuss this result, let us note  that the
numbers
$B_{k}(p)$  are represented by the sum over all possible  partitions
$\omega\in
\hat
\Omega_{2k}$ of the points
$1,2,
\dots, 2k$  to the set  of $k$ non-crossing pairs; each partition
is weighted by powers of $p$ (cf.(4))
$$
B_k(p) = \sum_{\omega\in \hat \Omega_{2k} } \prod_{(i,j)\in \o}p^{\vert
i-j\vert},
\eqno (8)
$$
such that 
$$
t_k = \sum_{\omega\in \hat \Omega_{2k} }1 =  {(2k)!\over k! (k+1)!} .
$$
Each non-crossing partition into pairs $\o\in \hat\O_{2k}$
can be identified with one of the half-plane rooted trees $T_k$
of $k$ edges \cite{S}. Thus $t_k$ represents the total number 
of the elements $\v T_k\v$.

Regarding definition (8), one can easily deduce that
$B_k(p)$ satisfy the following recurrent relations
$$
B_{k+1}(p) = \sum_{i=1}^{k+1} p^{2i-1} B_{i-1}(p) B_{k+1-i}(p)
\eqno (9)
$$
with the conditions $B_0(p)=B_1(p)=1$.
Indeed, following the reasoning by Wigner \cite{W},
let us consider the subsum of (8) over those partitions, where
the first point $1$ is paired with the last point $2k$;
we denote by $B'_k(p)$ the corresponding weighted sum. 
Then it is easy to observe that 
$$B'_k(p)=p^{2k-1} B_{k-1}(p).$$
To complete the derivation of (9), it remains to consider
the sum over partitions where $1$ is paired with $2i$.
Relations (8) and (9) answer the question of C. Mazza \cite{CM}
and theorem 1.1 establishes relations of (8) with
the product of random matrices of infinite dimensions.

In this connection, let us make one more remark about  relations
of our results with the non-commutative probability theory \cite{V}.  In
frameworks of this approach, the set of random matrices $A^{(r)}$ with $V(r)=
\d_{0,r}$ represents in the limit
$N\to\infty$  a family of free random variables 
$X^{(r)}$ with respect to the
mathematical expectation
$ \E \la \cdot\ra$ (6).
Free random variables represent a non-commutative analogue of
jointly independent scalar random variables.
In particular, 
according to the 
rules adopted to compute  the moments of 
free random variables (see e.g. \cite{NS}),
we recover the (ordinary) Catalan numbers
 $$
\lim_{N\to\infty}\la [A^{(1)}_N]^{2k}\ra = B_k(1) = t_k
$$
known since the pioneering work by Wigner 
on the eigenvalue distribuion of large random matrices \cite{W}. In  \cite{SP}
one can find a detailed analysis of the multiplicative functions
over the non-crossing partitions, where, in particular,  several
generalizations of the Catalan numbers appear. However, the q-Catalan numbers
are not present in
\cite{SP}.

In this context, conditions (1) can be regarded as a starting point
to define the non-commutative analogs $Y^{(r)}$ of the 
correlated scalar (gaussian) random variables.
The standard rules to compute the average value
$
\la Y^{(1)} \cdots Y^{(2k)}\ra
$ 
could be added there by conditions
$$
\la Y^{(s)}Y^{(t)}\ra = V(s-t).
$$

Returning to the generalization of the Catalan numbers (8),
let us note that as well as for the scalar case (5), 
 there
should be a  critical value  $p_c'$ in the sense that the limit $N\to\infty$
of (5) exhibits different behaviour as $k\to\infty$ in dependence whether 
$p>p_c'$ or $p<p_c'$. This can be explained by the observation that
the trees of $T_k$ that have a vertex of large degree are relatively
rare (see e.g. \cite{K}). However, the weight (8) with $p<1$ ascribes 
to such trees the probability greater than, say, to binary trees.
This makes the subject of  weighted non-crossing pairings reach and
interesting.

\vskip 1cm
{\it Proof of Theorem 1.1}. 

Our goal is to derive recurrent relation (9).
We rewrite (6) in the form
$$
 \la Q_{2k}^{(N)}(p) \ra  = {1\over N} \sum_{\{i\}}
A^{(1)}_{i_1i_2} A^{(2)}_{i_2i_3} \cdots A^{(2k)}_{i_{2k}i_1}\equiv
 \la Q_{(1, \dots, 2k)}^{(N)}(p) \ra 
$$
and compute the mathematical expectation with the help of (3) and then (1):
$$
\E \la Q_{2k}^{(N)}(p) \ra = 
{1\over N^{2}}  \sum_{\{i\}}
\sum_{l=2}^{2k} V(l-1) \times 
$$
$$
\E \left\{A^{(2)}_{i_2i_3} \cdots
A^{(l-1)}_{i_{l-1}i_l} 
\big( \d_{i_1i_l}\d_{i_2i_{l+1}} + \d_{i_1i_{l+1}}\d_{i_2i_{l}} \big)
A^{(l+1)}_{i_{l+1}i_{l+2}}\cdots
A^{(2k)}_{i_{2k}i_1}\right\}.
$$
In this relations we mean that for $l=2$ and $l=2k$ the
expression in curly brackets takes the forms
$A^{(3)} \cdots A^{(2k)}$ and $A^{(2)}\cdots A^{(2k-2)}$, respectively.

Then we obtain relation
$$
\E \la Q_{2k}^{(N)}(p) \ra = \sum_{l=2}^{2k} 
V(l-1) \E \la Q_{(2,3,\dots,l-1)}^{(N)}(p) \ra \la
Q_{(l+1,\dots,2k)}^{(N)}(p)\ra+
$$
$$
{1\over N} \sum_{l=1}^{2k} 
V(l-1) \E \la Q_{(2,\dots,l-1)}^{(N)}(p)
Q_{(2k, \dots, l+1)}^{(N)}(p)\ra.
\eqno (10)
$$
To complete the proof of Theorem 1.1, it remains to prove  
the following three 
items:

1) 
the moments of the normalized traces factorize
$$
\E \la Q_{(2,3,\dots,l-1)}^{(N)}(p) \ra \la
Q_{(l+1,\dots,2k)}^{(N)}(p)\ra  -
\E \la Q_{(2,3,\dots,l-1)}^{(N)}(p) \ra 
\ \E \la
Q_{(l+1,\dots,2k)}^{(N)}(p)\ra = o(1)
\eqno (11)
$$
in the limit $N\to\infty$;

2)  the odd moments are zero 
$$
\E \la Q_{(1,2,\dots,2m+1)}^{(N)}(p) \ra =0
\eqno (12)
$$
 and the average is invariant 
with respect to  simultaneous shifts of
all values  of the subscripts 
$$
\E \la Q_{(2m+1,2m+2,\dots,2k)}^{(N)}(p) \ra = 
\E \la Q_{(1,2,\dots,2k-2m)}^{(N)}(p) \ra
\eqno (13) 
$$
and finally 

3) the value of 
$$
\E \la Q_{(2,\dots,l-1)}^{(N)}(p)
Q_{(2k, \dots, l+1)}^{(N)}(p)\ra 
\eqno (14)
$$
 remains bounded as $N\to\infty$.

With these items in mind, we can easily derive from (10) that (7) takes place
and $B_k(p)$ satisfy (9).

Relations (12) and (13) trivially follow from the
definitions. Regarding (14), we observe that it is equal to  
$$
\E \la Q_{(1,\dots,l-2, 2k-1 \dots,l)}^{(N)}(p)\ra,
$$
where the number of factors $A$ is $2k-2$. Now
the question of its behavior as $N\to\infty$
is reduced to the problem of computing the expected value of (6)
where $2k-2$ factors are subjected to certain permutation.
 Thus, (14) can be  estimated in terms of
$\E \la Q_{2k-2}(p)\ra$ that is sufficient for us.

Let us describe the scheme of the proof of (11) based on  the following
standard procedure  (see e.g. \cite{BK}).
Let us denote $\xi^\circ = \xi- \E \xi$. Then 
we can write relations
$$
R^{(N)}_2([l+1,\dots,2k],[2,3,\dots,l-1]) \equiv
\E  \la
Q_{(l+1,\dots,2k)}^{(N)}(p)\ra^\circ 
\la Q_{(2,3,\dots,l-1)}^{(N)}(p) \ra=
$$
$$
\E \left\{\la
Q_{(l+1,\dots,2k)}^{(N)}(p)\ra^\circ {1\over N} \sum_{\{i\}}
A^{(2)}_{i_2i_3}  \cdots A^{(l-1)}_{i_{l-1}i_2}\right\}.
$$
To compute the last mathematical expectation, we use again the
identity (3) with $\gamma_l = A^{(2)}_{i_2i_3}$.
Repeating the same computations as above, we derive with the help of (1)
equality
$$
\E  \la
Q_{(l+1,\dots,2k)}^{(N)}(p)\ra^\circ 
\la Q_{2,3,\dots,l-1}^{(N)}(p) \ra= 
$$
$$
\sum_{j=3}^{l-1} V(2-j) \E  \left\{\la
Q_{(l+1,\dots,2k)}^{(N)}(p)\ra^\circ 
\la Q_{(3,\dots,j-1)}^{(N)}(p) \ra
\la Q_{(j+1,\dots,l-1)}^{(N)}(p) \ra
\right\}+
$$
$$
{1\over N} \sum_{j=3}^{l-1} V(2-j) \E  \left\{\la
Q_{(l+1,\dots,2k)}^{(N)}(p)\ra^\circ 
\la Q_{(3,\dots,j-1)}^{(N)}(p)  Q_{(j+1,\dots,l-1)}^{(N)}(p) \ra
\right\}+
$$
$$
{1\over N^2} \sum_{j= l+2}^{2k} V(2-j) 
 \E  \left\{\la
Q_{(2k,2k-1,\dots,j+1)}^{(N)}(p)   Q_{(3,\dots,l-1)}^{(N)}(p) Q_{(l+1,\dots,
j-1) }^{(N)}(p) 
 \ra
\right\}+
$$
$$
{1\over N^2} \sum_{j= l+2}^{2k} V(2-j) 
 \E  \left\{\la
Q_{(2k,2k-1,\dots,j+1)}^{(N)}(p)   Q_{(l-1,l-2,\dots,3s)}^{(N)}(p)
Q_{(l+1,\dots, j-1) }^{(N)}(p) 
 \ra
\right\}.
\eqno (15)
$$
This somewhat cumbersome relation has rather simple structure.
Indeed, using identity 
$$
\E \xi^\circ \zeta \mu = 
\E \xi^\circ \zeta \E \mu +
\E \xi^\circ \mu \E\zeta +
\E \xi^\circ \zeta^\circ \mu^\circ, 
$$
we observe that $R_2= \E Q^\circ Q^\circ$ is expressed as the sum of terms
of the following forms: 
$R_2 \E Q$,
$R_3 = \E Q^\circ Q^\circ Q^\circ$, $R_2/N$ and $\E Q/N^2$.
Also we see that in (15) the total number of factors $A$ is decreased by 2.

Thus, we conclude that repeating this procedure with respect to all terms
$R_m$ in (15), we obtain the finite number (depending on $k$) of terms that 
involve  products of $V$ and $1/N$.
The only expressions that has no factor $1/N$ are 
$$
R_{m} = \E \left\{\prod_{j=1}^{m} \la
Q_{(\alpha_j,\beta_j)}^{(N)}\ra^\circ\right\}, m\le k/3.
$$
This averaged product can be treated as before with the help of (3).
We obtain relations that express $R_m$ in terms of the sum
of $R_{m-1}$ and terms that have factors $1/N$. Taking into account that $k$
is finite, we arrive at (11). 
Theorem 1.1 is proved.
\hfill $\Box$


\begin{thebibliography}{99}

\bibitem{BK} A. Boutet de Monvel and A. Khorunzhy,
On the norm and eigenvalue distribution of large random matrices,
{\it Ann. Probab.} {\bf 27} (1999) 913-944

\bibitem{K} A. Khorunzhy, Sparse random matrices: spectral edge and
statistics of rooted trees, {\it Adv. Appl. Probab.} {\bf 33} (2001)
124-140


\bibitem{CM} C. Mazza, {\it private communication}

\bibitem{M} M L Mehta, {\it Random Matrices}, Acad. Press (1991)


\bibitem{SP} R. Speicher, Multiplicative functions on the lattice
of non-crossing partitions and free convolution,
{\it Math. Annalen} {\bf 298} (1994) 611-628

\bibitem{NS} P. Neu and R. Speicher, Rigorous mean field model for
CPA: Anderson model with free random variables.
{\it J. Stat. Phys.} {\bf 80} (1995) 1279-1308

\bibitem{S} R.P. Stanley, {\it Enumerative Combinatorics}, vol.II,
Cambridge University Press (1999)



\bibitem{V} D.V.  Voiculescu, K. Dykema and A. Nica,
{\it Free Random Variables}, CRM Monograph Series,
No. 1, Provodence, RI (1992)

\bibitem{W} E. Wigner, Characteristic vectors of bordered matrices with
infinite dimensions, {\it Ann. of Math.} {\bf 62} (1955) 548-564 

\end{thebibliography}
\end{document}